          \documentclass[11pt]{amsart} \usepackage{latexsym,amsmath,amsthm,amssymb} 
                                       \usepackage{graphicx}
               \def\bb{\Bbb}      \def\bl{\bar{l}} 
                   \def\bm{\bar{\mu}} 
                   \def\bn{\bar{n}}
               \def\bQ{\bb{Q}}    \def\bz{\bar{0}}
               \def\bZ{\bb{Z}}    

               \def\emp{\emptyset}      \def\prd{\times}

               \def\rt{\sqrt{t}}       \def\wtd{\widetilde} 

                 \def\ptl{\partial}       

               \def\Cup{\bigcup}         \def\cmpl{S^3 \setminus}  

               \def\ccc{\cup \cdots \cup}

               \def\ti{t^{-1}}    \def\rtt{(\rt-\rti)}
                  \def\rti{\sqrt{t}^{-1}}

               \def\cd{\cdot}     \def\cds{\cdots}  
                      
                    \def\vds{\vdots}

               \def\lk{\mathrm{lk}} 
               \def\det{\mathrm{det}} 

               \def\CwyK{\na_K(z)}  \def\CwyKS{\na_{K_S}(z)}   \def\AlxK{\Del_K}
               \def\CwyL{\na_L(z)}  \def\CwyLS{\na_{L_S}(z)}   \def\AlxL{\Del_L}

               \def\sumkzi{\sum\limits_{k=0}^{\infty}}
               \def\sumkzn{\sum\limits_{k=0}^{n}}
               \def\sumkai{\sum\limits_{k=1}^{\infty}}
               \def\sumkam{\sum\limits_{k=1}^m}
               \def\sumnzi{\sum\limits_{n=0}^{\infty}}

               \def\iab{(i=1,2)}
               
               \def\hZ{{h \in \bZ}}
               \def\kZ{{k \in \bZ}} 

               \def\tF{\wtd{F}}                    \def\mJi{\mJ_i}
               \def\tJ{\wtd{J}}   \def\JF{\tJ^F}   \def\mJF{\wtd{\mJ}^F}
                                   \def\mJiF{\wtd{\mJ}_i^F}

               \def\Ja{J_1}       \def\tJa{{\tJ_1}}   \def\JaF{{\tJa^F}} 
               \def\Jb{J_2}       \def\tJb{{\tJ_2}}   \def\JbF{{\tJb^F}} 
               \def\Ji{J_i}       \def\tJi{{\tJ_i}}   \def\JiF{{\tJi^F}} 
               \def\Jj{J_j}       \def\tJj{{\tJ_j}}   \def\JjF{{\tJj^F}} 
               \def\Jk{J_k}       \def\Jip{{\Ji^+}}   \def\JipF{{\tJi^{+F}}}
               \def\Jm{J_m}

               \def\KJ{K       \cup J} 
               \def\JaJb{\Ja   \cup \Jb}
               \def\KJaJb{K    \cup \JaJb}
               \def\KzJaJm{K_0 \cup \Ja \ccc \Jm}

            \def\fD{D_F}      
            \def\pD{D_{F_+}}    \def\pnD{D_{F_\pm}} 
            \def\nD{D_{F_-}}    

               \def\pJia{{J_i^{(  1)+}}}   \def\tpJa{\tJa^{(1)+}} 
                                           \def\pJak{\Ja^{(k)+}}
               \def\pJim{{J_i^{(n-1)+}}}  
               \def\pJin{{J_i^{(n  )+}}}

               \def\nJin{{J_i^{(\bar{n  })-}}} 
             
            \def\nJbl{\Jb^{(\bar{l})-}}

            \def\MO{M_\Om}        \def\MOt{M_\Om^T}       
                
            \def\VJa{V_{\Ja}}     \def\VJi{V_{\Ji}}
            \def\VJb{V_{\Jb}}     \def\VJj{V_{\Jj}}    
            \def\VJbt{V_{\Jb}^T}  \def\VJm{V_{\Jm}}
         
            \def\pa{p_1} \def\pap{p_1^+} \def\qa{q_1} \def\qap{q_1^+} 
            \def\pi{p_i}  \def\qi{q_i}  
            \def\pj{p_j} \def\pjp{p_j^+} \def\qj{q_j} \def\qjp{q_j^+} 
            \def\pg{p_g} \def\pgp{p_g^+} \def\qg{q_g} \def\qgp{q_g^+} 

            \def\VRV{V R^{-1}(t) V^T} 
            \def\pq{\{ \pa, \ldots, p_g, q_1, \ldots, q_g \}} 
            \def\pqr{\{\pa, \ldots, p_g,  q_1, \ldots, q_g, r_1, \ldots, r_m \}}

            \def\alfz{\al_F^0}   \def\akab{\al_F^k(J_1,J_2)}       
            \def\alfa{\al_F^1}   \def\akij{\al_F^k(J_i,J_j)}
            \def\alfk{\al_F^k}   \def\anab{\al_F^n(J_1,J_2)}       
            \def\alfn{\al_F^n}   
                                 \def\azab{\al_F^0(J_1,J_2)}

            \def\beaz{\be_F^{1,0}} 
                  
            \def\bekl{\be_F^{k,l}}      \def\bklab{\be_F^{k,l}(\Ja,\Jb)} 
            \def\benz{\be_F^{n,0}}
          
            \def\lkx{\lk_X}
            \def\lkxab{\lkx(\tJa,\tJb)}    
            \def\lkxabf{\lkx(\JaF,\JbF)}  
            \def\lkxabd{\lkx(\tJa^{F'},\tJb^{F'})}  

\def\jajb{\lk(J_1,J_2)} 
 \def\jipa{\lk(J_i,p_1)}  
 \def\jipj{\lk(J_i,p_j)} 
 \def\jipg{\lk(J_i,p_g)}  
 \def\jiqa{\lk(J_i,q_1)}  
 \def\jiqj{\lk(J_i,q_j)}
 \def\jiqg{\lk(J_i,q_g)}

\def\papa{\lk(p_1^+,p_1)} \def\papg{\lk(p_1^+,p_g)} 
\def\paqa{\lk(p_1^+,q_1)} \def\paqg{\lk(p_1^+,q_g)}
\def\pgpa{\lk(p_g^+,p_1)} \def\pgpg{\lk(p_g^+,p_g)} 
\def\pgqa{\lk(p_g^+,q_1)} \def\pgqg{\lk(p_g^+,q_g)}
\def\qapa{\lk(q_1^+,p_1)} \def\qapg{\lk(q_1^+,p_g)} 
\def\qaqa{\lk(q_1^+,q_1)} \def\qaqg{\lk(q_1^+,q_g)}
\def\qgpa{\lk(q_g^+,p_1)} \def\qgpg{\lk(q_g^+,p_g)} 
\def\qgqa{\lk(q_g^+,q_1)} \def\qgqg{\lk(q_g^+,q_g)}

          \def\deta{\left| \begin{array}{c}}      \def\mtxa{\left( \begin{array}{c}} 
          \def\detb{\left| \begin{array}{cc}}     \def\mtxb{\left( \begin{array}{cc}} 
          \def\detc{\left| \begin{array}{ccc}}    \def\mtxc{\left( \begin{array}{ccc}} 
          \def\detd{\left| \begin{array}{cccc}}   \def\mtxd{\left( \begin{array}{cccc}} 
          \def\dete{\left| \begin{array}{ccccc}}  \def\mtxe{\left( \begin{array}{ccccc}}
          \def\detf{\left| \begin{array}{cccccc}} \def\mtxf{\left( \begin{array}{cccccc}} 
          \def\enddet{\end{array}  \right|}       \def\endmtx{\end{array}  \right)} 

          \def\bgp{\begin{proof}} 

                              \def\al{\alpha}
                              \def\be{\beta} 
                                     \def\Gam{\Gamma}  
                              \def\de{\delta}       \def\Del{\Delta}

                              \def\et{\eta}

                              \def\la{\lambda}      
                              
                              \def\si{\sigma}        
                              \def\ta{\tau}         \def\tam{\ta^{-1}}
                                    
                              \def\ph{\phi}

                              \def\om{\omega}       \def\Om{\Omega}
                              \def\na{\nabla}

           \def\mJ{\mathcal{J}}

          \def\pvc{\par\vspace*{3mm}} \def\phc{\par\hspace{3cm}} 
           \def\phd{\par\hspace{4cm}}  
          \def\pve{\par\vspace*{5mm}}

          \def\phc{\par\vspace{3mm}\hspace{5mm}}
           \def\pvcn{\par\vspace{3mm}\noindent}
          \def\phcb{\par\vspace{3mm}\hspace{2cm}} \def\pvdn{\par\vspace{4mm}\noindent}
          \def\phcc{\par\vspace{3mm}\hspace{3cm}} 
          \def\phcd{\par\vspace{3mm}\hspace{4cm}} \def\phcdh{\par\vspace{3mm}\hspace{4.5cm}}

          \newtheorem{Thm}{Theorem}[section]  \newtheorem{Cor}[Thm]{Corollary}
          \newtheorem{Lem}[Thm]{Lemma}         
                  
            \newtheorem{Prop}[Thm]{Proposition}
               \newtheorem{Rem}[Thm]{Remark}

                  \title{A factorization of the Conway polynomial 
                                  and covering linkage invariants } 
                 \author{Tatsuya Tsukamoto \\
                         Department of Mathematical Sciences, 
                         School of Science and Engineering,  \\
                         Waseda University, 3-4-1 Okubo Shinjuku-ku, 
                         Tokyo 169-8555 JAPAN\\
\\
                         and\\
\\
                         Akira Yasuhara \\
                         Department of Mathematics,
                         Tokyo Gakugei University, \\
                         Nukuikita 4-1-1, Koganei, 
                         Tokyo 184-8501 JAPAN\\} 
                 \thanks{The first named author acknowledges partial support by 
                         JSPS Research Fellowships for Young Scientists.}
                 \date{}
\begin{document} 
                                   \maketitle \baselineskip=15pt 
                        
%
%
%

\begin{abstract} J.P. Levine showed that the Conway polynomial of a link is a product
                 of two factors: one is the Conway polynomial of a knot which is
                 obtained from the link by banding together the components; and
                 the other is determined by the $\bm$-invariants of a string link
                 with the link as its closure.
                 We give another description of the latter factor: the determinant
                 of a matrix whose entries are linking pairings in the infinite
                 cyclic covering space of the knot complement, which take values
                 in the quotient field of ${\Bbb Z}[t,t^{-1}]$.
                 In addition, we give a relation between the Taylor expansion of
                 a linking pairing around $t=1$ and derivation on links which is
                 invented by T.D. Cochran.
                 In fact, the coefficients of the powers of $t-1$ will be
                 the linking numbers of certain derived links in $S^3$.
                 Therefore, the first non-vanishing coefficient of the
                 Conway polynomial is determined by the linking numbers in $S^3$.
                 This generalizes a result of J. Hoste.
                 \end{abstract}

%
%
%

                                  \section{Statement of Results}
                                    \label{sec:Introduction}

               Let $L$ be an oriented $(m+1)$-component link ($m \geq 1$) 
               in the $3$-sphere $S^3$. Throughout the paper knots and links 
               are assumed to be oriented. It is interesting but difficult 
               in general to figure out the structure of a link invariant.
               The Conway polynomial $\CwyL$ of $L$ , the most popular 
               polynomial link invariant, was not an exception, either.
               However, J.P. Levine \cite{JL-CMH99} made a breakthrough.
               Namely he showed that there is a  following relationship 
               between $\CwyL$ and $\CwyK$, where $K$ is a knot obtained 
               from $L$ by banding together the components and $\Gam(z)$ 
               is a power series in $z$ which depends on the choice of bands: 
               $$\CwyL = \CwyK \ \Gam(z).$$ 
               Viewing the choice of bands as the choice of a string link 
               represenation of $L$, Levine descrived $\Gam(z)$ by
               $\bm$-invariants of the string link as follows. 

\begin{Thm}    $($\cite{JL-CMH99} 
               $\mathrm{Theorem} \ 1)$ \label{thm:FactorizationofLevin}   
               Let $S$ be a string link, with closure $L_S$ and knot closure 
               $K_S$. Then we have that $\CwyLS = \CwyKS \Gam_S(z)$, where 
               $\Gam_S(z)$ is a power series given by the formula:
               $$\Gam_S(z) = (u+1)^{e/2} \ \det(\la_{ij}(u)) \ \ \ with \
               z=u/\sqrt{u+1}, \ e= \left\{ \begin{array}{r@{\quad \quad}l} 
               0 & if \ m \ is \ even \\ 
               1 & if \ m \ is \ odd  \end{array} \right.and:$$ 
               $$\la_{ij}(u)=\sum_{r=0}^\infty (\sum_{i_1,\cds,i_r}
                 \bm_{i_1,\cds,i_r,j,i}(S)) \ u^{r+1} \ \ \ (1 \leq i,j \leq m).$$
               \end{Thm}

         \pvcn In this paper we give another description of $\Gam(z)$ by viewing 
               the choice of bands as the choice of a Seifert surface for $L$. 
               Namely $\Gam(z)/z^m$ is the determinant of a matrix whose entries 
               are linking pairings in the infinite cyclic covering space of 
               the complement of the knot obtained from $L$ by banding together 
               the components.
               To state the theorem first let us recall the definitions.

         \pvcn Let $F_0$ be a Seifert surface for $L$ and 
               $\Om$ a basis for $H_1(F_0)$.
               A {\it Seifert form} is a map $\si: H_1(F_0) \prd H_1(F_0) \to \bZ$ 
               defined by $\si(x,y)= \lk(x^+,y)$, where $x^+$ is $x$ 
               pushed into the positive normal direction of $F_0$.
               Let $\MO$ be a Seifert matrix representing $\si$ with respect to $\Om$.
               Define the {\it potential function} $\AlxL(t) = \det(\rt \MO - \rti \MOt)$. 
               It is known that $\AlxL(t)$ is a polynomial in $\rtt$. 
               Then we may define the {\it Conway polynomial} $\CwyL$ of $L$ by 
               $\na_L\rtt=$ $\AlxL(t)$.

         \pvcn Next let $K$ be a knot and $F$ a Seifert surface for $K$.
               Consider two knots $\Ja$ and $\Jb$ in $\cmpl F$.
               Let $\ph: X \to \cmpl K$ be the infinite cyclic cover and 
               let $\ta$ be a covering transformation that shifts $X$ by one 
               along the positive direction with respect to $K$. 
               Take a component $\tF$ of $\ph^{-1}(F)$ and 
               let $\JiF$ be the component of $\ph^{-1}(\Ji)$ 
               between $\tF$ and $\ta\tF$ ($i=1,2$).
               Then the {\it linking pairing} $\lkxabf$ of 
               $\JaF$ and $\JbF$ in $X$ is defined as follows, 
               where $H$ is a $2$-chain in $X$ such that 
               $\ptl H = \Cup_\kZ c_k \ta^k \JaF$:
               
               $$\lkxabf = \frac{\sum_\hZ (H \cd \ta^h \JbF) \ t^h}
                                {\sum_\kZ c_k \ t^k} \in \bQ(\bZ[t,t^{-1}]).$$

\begin{figure}[htbp!] 
               \centerline{\includegraphics[scale=.6, bb=98 285 496 674]
               {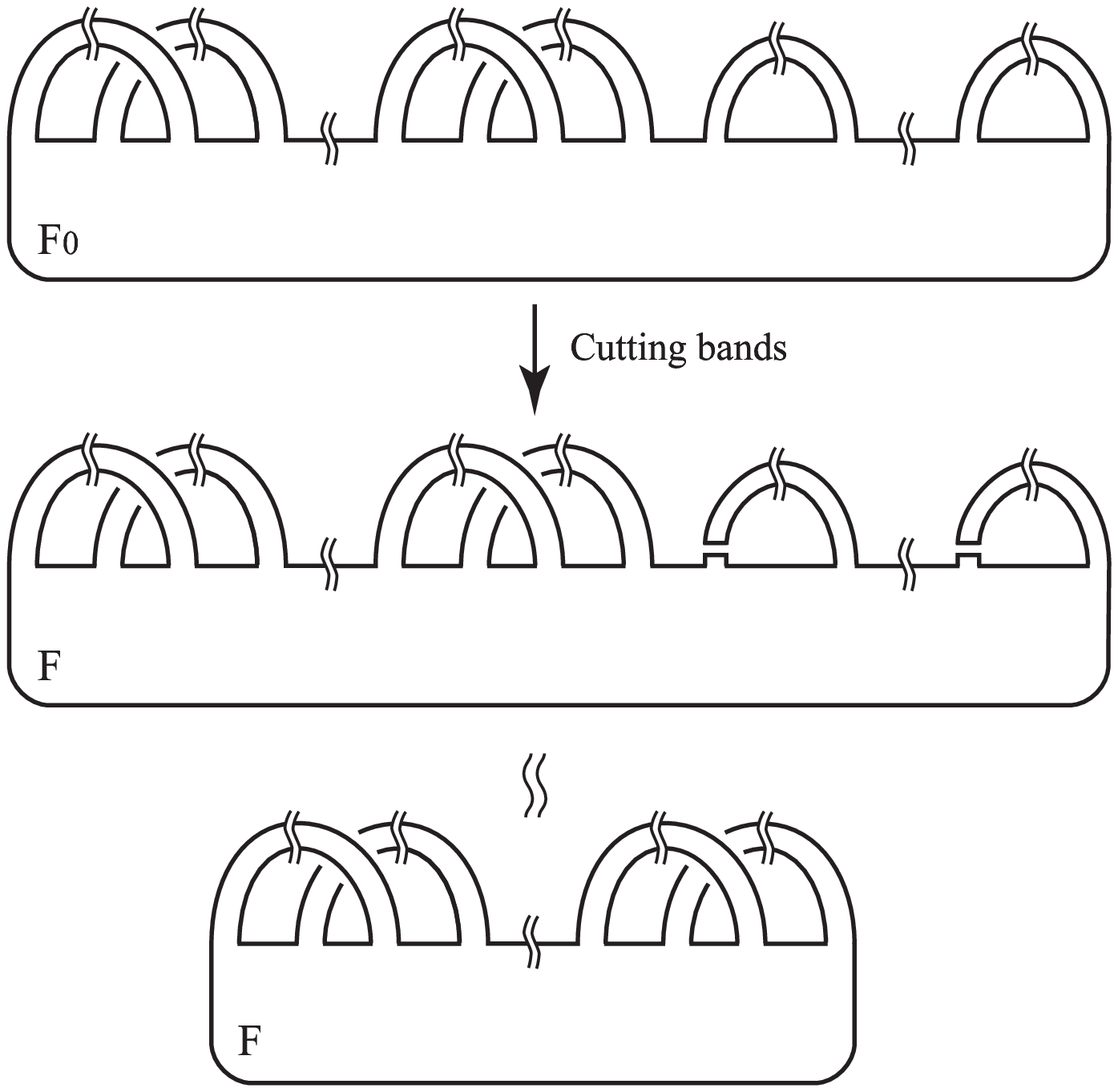}}\caption{}\label{fig:fusions}
               \end{figure}

         \pvdn Now consider an $(m+1)$-component link $L = \KzJaJm$. 
               Let $F_0$ be a Seifert surface for $L$ and 
               $\Om = \pqr$ be a basis for $H_1(F_0)$, 
               where $\pq$ is symplectic, i.e.
               $\pi \cd \pj=$ $\qi \cd \qj=0$ and $\pi \cd \qj = \de_{ij}$,
               and $r_i$ is represented by $\Ji$. 
               Let $F$ be a surface obtained from $F_0$ by cutting bands 
               which are presented by $r_1, \ldots, r_m$ 
               (see Figure \ref{fig:fusions}).  
               Then the boundary of $F$ is a knot, say $K$.
               Consider the infinite cyclic covering space $X$ of $\cmpl K$
               and define $\JiF$ similarly to the above.
               Then we have the following, where note that $\Jip$ is 
               $\Ji$ pushed into the positive normal direction of $F_0$.

\begin{Thm}    \label{thm:FctznOne} 
               $\CwyL  = z^m$ $\CwyK$ $\det(p_{ij})$, where
               \[p_{ij} = \lkx(\JipF,\JjF) = \left\{ \begin{array}{l l}
                         \lkx(\JipF,\JiF) & if \ i=j, \\
                         \lkx(\JiF, \JjF) & if \ i \neq j. \end{array} \right.\]
\end{Thm}
               \pve          

          \pve Next we take a look at entries of $(p_{ij})$, linking pairings.
               T.D. Cochran \cite{TC-CMH85} showed that the coefficients of
               the Maclaurin series of $\lkx(\JiF,\mJiF)$ are the linking number
               of $K$ and a ``derivative" of $\Ji$, where $\mJi$ is a parallel 
               copy of $\Ji$ with $\lk(\Ji,\mJi)=0$. 
               Derivation is defined also in \cite{TC-CMH85} for 
               a $2$-component link with its linking number vanished.
               Here following and slightly modifying his idea, we show 
               a similar theorem for a Taylor series of $\lkx(\JipF,\JjF)$ 
               around $t=1$. First let us define derivation on a link 
               with respect to a Seifert surface for a knot
              (We will apply derivation to a knot. However since the
               derivative of the knot might be a link, here we define
               derivation for a link).
               
\begin{figure}[htbp!] 
               \centerline{\includegraphics[scale=.6, bb=143 503 433 721]
               {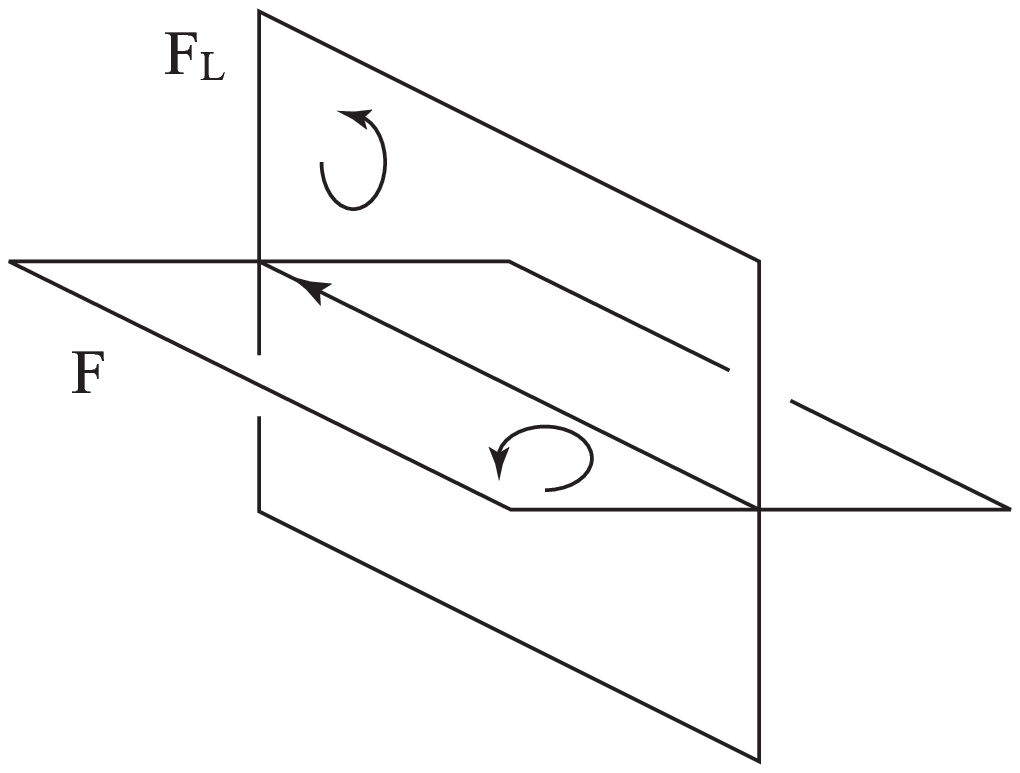}}\caption{}\label{fig:orientation}
               \end{figure}

         \pvcn Let $K$ be a knot and $F$ a Seifert surface for $K$.
               Let $L$ be a link in $\cmpl F$ and $F_L$ a Seifert surface 
               for $L$ in $\cmpl K$. 
               We may assume that the non-empty intersection of 
               $F_L$ and $F$ consists of simple closed curves on $F$. 
               Then from the following proposition we can define 
               {\it derivation} $\fD(L)$ of $L$ with respect to $F$ 
               as the homology class of the intersection 
               $F_L \cap F$ on $F$ with the orientation following the order 
               $F_L$, $F$ (see Figure \ref{fig:orientation}).
               If the intersection is empty, then we define $\fD(L)=0$.

\begin{Prop}   \label{prop:derivative}     
               Let $K$ be a knot $($resp. $L$ be a link$)$ in $S^3$ such 
               that there is a Seifert surface $F$ for $K$ in $\cmpl L$ 
               $($resp. $F_L$ for $L$ in $\cmpl K)$. Then the homology class 
               of the intersection of $F_L$ and $F$ on $F$ is determined 
               by the homology class of $L$ in $\cmpl F$. 
               \end{Prop}
          \bgp 
               Take a symplectic basis $\pq$ for $H_1(F)$. Let $c$ be the 
               intersection of $F_L$ and $F$. Since $c$ is on $F$ and the
               basis is symplectic, $c$ can be wrote as               
               $c = \sum ((c \cd \qi) \pi - (c \cd \pi) \qi)$. Now since 
               $\pi$ and $\qi$ are on $F$, 
               $c \cd \pi=    F_L \cd \pi = \lk(L, \pi)$ and 
               $c \cd \qi=    F_L \cd \qi = \lk(L, \qi)$. Thus we obtain that 
               $c = \sum (\lk(L, \qi) \pi - \lk(L, \pi) \qi)$, which implies that 
               $c$ is determined only by   $\lk(L, \pi)$ and $\lk(L, \qi)$.  
\end{proof}
               
         \pvcn Then we define {\it positive derivation} $\pD^n$
                       (resp. {\it negative derivation} $\nD^n$) on $L$
               with respect to $F$ for a natural number $n$ as follows:
               let $\pnD^1(L)=\pnD(L)=\fD(L)$ and for $n \geq 2$
               push the $n-1$ times derived link into the positive 
               (resp. negative) normal direction of $F$ and apply $\fD$, i.e.
               $\pnD^n(L)=\fD(\pnD^{n-1}(L)^\pm)$.
               This is well-defined also from Proposition \ref{prop:derivative}.  
               Here note that there is a Seifert surface for $\pnD^{n-1}(L)^\pm$ 
               in $\cmpl K$, since $\pnD^{n-1}(L)$ is on $F$ and thus 
               its each component has the vanished linking number with $K$. 
               We also denote $\pD^n(L)$ by $L^{F(n)}$ (resp. $\nD^n(L)$ 
               by $L^{F(\bn)}$) and we omit $F$ if no confusion is expected.
               \pvc

          \pvc Now let $\KJaJb$ be a $3$-component link with $\lk(K, \Ji)=0$ 
               $(i=1,2)$ and $F$ a Seifert surface for $K$ in $\cmpl (\JaJb)$.
               For non-negative integers $k$ and $l$, we define $\bekl$ 
               by $\bekl(\Ja,\Jb)= \lk(\Ja^{F(k)+}, \Jb^{F(\bl)-})$, where 
               $\Ja^{F(0)}=\Ja$ and $\Jb^{F(\bz)} =\Jb$ and 
               $\lk(0,*)=$ $\lk(*,0)=0$. The function $\bekl$ is 
               well-defined also from Proposition \ref{prop:derivative}. 
               Since we will see that $\bekl(\Ja,\Jb)=(-1)^l \be_F^{k+l,0}(\Ja,\Jb)$
              (Corollary \ref{cor:beta}), we define $\alfn$ by
               $\anab = \benz(\Ja,\Jb)$ for $n \geq 0$.
               Then we can obtain the following theorem by the arguments
               similar to that in the proof of Theorem 7.1 in \cite{TC-CMH85}.

\begin{Thm}    \label{thm:ICC} $(c.f.$ \cite[Theorem 7.1]{TC-CMH85} $)$
               The Taylor expansion of the linking pairing of $\JaF$ and $\JbF$ 
               around $t=1$ is descrived as follows:
               $$\lkxabf = \sumkzi (-1)^k \akab (t-1)^k.$$
               \end{Thm}
               \pvcn
               In general, $\anab$ depends on the choice of $F$. 
               However, the first non-vanishing one does not. 
               Namely we have the following as a corollary of the above theorem.

\begin{Cor}    \label{cor:Fnva} If $\akab$ vanishes for all $k < n$, 
               then $\anab$ is an invariant of $\KJaJb$. 
               \end{Cor}

\begin{Rem}   {\em Note that 
               $\anab=\lk(\Ja^{F(n)+},\Jb)= \lk(F_{\Ja^{F(n-1)+}}\cap F,\Jb)$
               $(n\geq 1)$, which is the algebraic sum of the triple points 
               of $F_{\Ja^{F(n)+}}\cap F\cap F_{\Jb}$. 
               Thus $\anab$ is not only the linking number for the derived link 
               $\Ja^{F(n)+} \cup \Jb$ but also the triple Milnor $\bm$-invariant
               $\bm(\Ja^{F(n-1)+},K,\Jb)$
               for the delived link $\Ja^{F(n-1)+}\cup K\cup \Jb$.}
               \end{Rem}

\begin{Rem}    {\em We can also define another function 
               $\lk(\Ja^{F(k)+}, \Jb^{F(l)})$.
               In fact, for a $2$-component link $\KJ$ with $\lk(K,J)=0$,
               $\lk(J^{F(k)+}, J^{F(k)})$ is equivalent to the function 
               defined by Cochran in \cite{TC-CMH85}. Here note that his 
               derivation is defined on the weak-cobordism class of $\KJ$
               and then the derivation is independent from the choice of a 
               Seifert surface for $K$.
               If $\mJ$ is parallel to $J$ and $\lk(J,\mJ)=0$, then
               the pairing $\lkx(\JF,\mJF)$ does not depent on the choice
               of $F$ and thus it is an invariant of $\KJ$. It is called
               the {\it Kojima-Yamasaki's} $\et$-{\it function} 
               $\et(K,J:t)$ $($\cite{KY-IM79}$)$. Cochran gave an expansion
               of this function: $\lkx(\JF, \mJF) = 
               \sum_{k=1}^\infty \lk(J^{F(k)+}, J^{F(k)}) x^k$, 
               where $x=(t-1)(\ti-1)$.}
               \end{Rem}

          \pvc Now again consider the Conway polynomial $\CwyL$ of $L=\KzJaJm$
               and take a Seifert surface $F_0$ for $L$ as before.
               It is known that the $i$-th coefficient 
               of $\CwyL$ ($i\leq m-1$) vanishes.
               Then J.Hoste \cite{JH-PAMS85} (and F.Hosokawa 
               \cite{FH-OMJ58} for the absolute value) showed that 
               the coefficient of $z^m$ is $\det(\lk(\Jip,\Jj))$.

\noindent      If $L$ be algebraically split, then Levin \cite{JL-Knots96} 
               showed that moreover the coefficients of $z^m$, $\cds$, $z^{2m-3}$
               vanish and that the coefficient of $z^{2(m-1)}$ is
               $\det(\bm(\Ji, \Jj, K_0) + \sum_{k \neq i} \bm(\Ji,\Jj,\Jk))$.

        \pvcn  Here we know from Theorem \ref{thm:FctznOne} and Theorem \ref{thm:ICC}
               that the first non-vanishing coefficient of $\CwyL$ are determined 
               by $\akij$, since $\na_K(0)=1$. Moreover under a certain condition, 
               we can obtain the explicit formulae as in the following corollary. 
               Since $\alfz(\Jip,\Jj)=\lk(\Jip,\Jj)$, the corollary generalizes 
               the above result of Hoste. The assumption in the case
               $n=1$ implies that $L$ is algebraically split.
               Therefore the corollary also generalizes the above result of Levine,
               since we have that
               $$\alfa(\Jip,\Jj)=\lk(F_{\Jip} \cap F, \Jj)
                                =\bm(\Jip,K,\Jj)
                                =\bm(\Jip,K_0,\Jj)+\sumkam \bm(\Jip,J_k,\Jj).$$

\begin{Cor}    \label{cor:Fnvc} 
               Assume that $n=0$ or that $n \geq 1$ and $\alfk(\Jip,\Jj)$ vanishes 
               for any $k < n$ and any pair of $i$ and $j$ $(1 \leq i \leq j \leq m)$. 
               Then $\CwyL$ is divisible by $z^{m(n+1)}$ and 
               the coefficient of $z^{m(n+1)}$ is $\det((-1)^n \alfn(\Jip,\Jj))$. 
               \begin{flushright}$\Box$\end{flushright}\end{Cor}

\noindent      The paper is organized as follows. In Section \ref{sec:Derivative}
               we study derivations on links and linking pairings.
               Theorem \ref{thm:FctznOne} is proved in Section \ref{sec:Conway}.
               In the last section, we present a calculation of $\bekl$ and 
               the Taylor expansion of the inverse of the Alexander matrix 
               $tM-M^T$ by a seifert matrix.

%
%
%
%
%
%
%
%
%
%
 
                \section{Derivation on links and Covering linkage invariants}
                \label{sec:Derivative}

%
%
%
%
%
%
%
%
%

        \pvc   In this section we prove Theorem \ref{thm:ICC} and 
               Corollary \ref{cor:Fnva}. To do so, we need the following Lemma, 
               which is clear from the definition of the linking pairing.

\begin{Lem}    \label{lem:ICC}
               Let $K$ be a knot and let $\ph: X \to \cmpl K$ be the infinite 
               cyclic cover with a covering transformation $\ta$ that shifts $X$ 
               by one along the positive direction with respect to $K$. 
               Let $L$, $L_1$, and $L_2$ be links in $X$. 
               We assume that there is a $2$-chain $H$ such that 
               $\ptl H = L_1 - L_2$ and $L \cap \ta^l H = \emp$ 
               for any $l$ $(\neq 0)$. Then we have the following. \begin{enumerate}
\item[$(i)$] $\lkx(L_1,L) = \lkx(L_2,L)+ H \cd L$, 
\item[$(ii)$] $\lkx(\ta L_1,L_2)=$ $\lkx(L_1, \tam L_2)=$ $t \lkx(L_1,L_2)$.
 \end{enumerate}
               \begin{flushright}$\Box$\end{flushright}\end{Lem}

         \pvcn In the following proofs, we omit $F$ of $\JiF$ and
               write $\akab$ as $\alfk$ in short for a convenience
               unless otherwise stated.. 
\pvcn
{\it Proof of Theorem} \ref{thm:ICC}.
               Let $H_\tJa$ be the component of 
               $\ph^{-1}(F_{\Ja})-\bigcup \ta^k \tpJa$ with
               $\ptl H_\tJa = \tJa \cup (-\tpJa \cup \ta \tpJa)$,
               where $F_{\Ja}$ is a Seifert surface for $\Ja$ in $\cmpl K$. 
               Then we have that 
               $H_\tJa \cap \ta^l \tJb = \emp$ for any $l$ $(\neq 0)$, and that  
               $H_\tJa \cd \tJb  
               =F_{\Ja} \cd \Jb = \jajb$.
               Thus from Lemma \ref{lem:ICC} (i) we obtain the following:
               $$\lkxab = \lkx(\tpJa - \ta \tpJa, \tJb) + \lk (\Ja,\Jb)
               =(1-t) \ \lkx(\tpJa, \tJb) + \alfz.$$
         \pvcn It is inductively and similarly shown that, 
               for any natural number $n$
         \pvcn $$\lkxab = (-1)^{n+1} (t-1)^{n+1} \lkx(\tJa^{(n+1)+}, \tJb) 
               + \sumkzn (-1)^k \alfk (t-1)^k$$ 
         \pvcn Then we have that
               $$\frac{d^n(\lkxab)}{dt^n} |_{t=1} = (-1)^n n! \alfn.$$
                \begin{flushright}$\Box$\end{flushright}
\pvcn
{\it Proof of Corollary} \ref{cor:Fnva}.
               The statement holds for $n=0$, since $\azab= \jajb$ by definition. 
               Now assume that $n \neq 0$. From the induction hypothesis, 
               we have that $\al_{F'}^l=0$ for another Seifert surface $F'$ 
               for $K$ in $\cmpl \JaJb$ and any $l$ ($< n$). 
               Thus we have the following from Theorem \ref{thm:ICC}.
               $$\lkxabf = (-1)^n     \al_F^n        (t-1)^n        
                         + (-1)^{n+1} \al_F^{n+1}    (t-1)^{n+1} + \cds  $$ 
           and $$\lkxabd = (-1)^n     \al_{F'}^n     (t-1)^n     
                         + (-1)^{n+1} \al_{F'}^{n+1} (t-1)^{n+1} + \cds. $$
         \pvcn From Lemma \ref{lem:ICC} (ii), 
               we have that $\lkxabd = t^m \lkxabf$ 
               for some $m$. Hence we have that:  
               $$ (-1)^n  n! \al_{F'}^n = \frac{d^n(t^m \lkxabd)}{dt^n} |_{t=1}
                                        = \frac{d^n(    \lkxabf)}{dt^n} |_{t=1}
                = (-1)^n  n! \al_F^n.$$
                                       \begin{flushright}$\Box$\end{flushright}

%
%
%
%
%
%
%

                        \pve \section{A Proof of Theorem \ref{thm:FctznOne}}
                                       \label{sec:Conway}
%
%
%
%
%
%
%
   
\noindent     In this section we give a proof of Theorem \ref{thm:FctznOne}.
\pvcn
{\it Proof of Theorem \ref{thm:FctznOne}.}
               Let $A= (\si(r_i,r_j))$ and define $M$, $\VJi$, 
               and $V$ as follows: 
               $$M = \mtxf  
               \papa & \cds & \papg & \paqa & \cds & \paqg  \\
               \vds  &      & \vds  & \vds  &      & \vds   \\
               \pgpa & \cds & \pgpg & \pgqa & \cds & \pgqg  \\
               \qapa & \cds & \qapg & \qaqa & \cds & \qaqg  \\
               \vds  &      & \vds  & \vds  &      & \vds   \\
               \qgpa & \cds & \qgpg & \qgqa & \cds & \qgqg  \\ \endmtx,$$ 
        \pvcn
               $$\VJi   =(\jipa, \cds, \jipg, 
                          \jiqa, \cds, \jiqg), \ \ \mathrm{and} \ \ 
                     V  = \mtxa \VJa \\ \vds \\ \VJm \\ \endmtx.$$
        \pvcn  
               Then we have that $\MO= \mtxb M & V^T \\ V & A \\ \endmtx$
               and using the notation we can calculate potential 
               function $\AlxL(t)$, where we let $R(t) = \rt M - \rti M^T$. 
               Note that $\det R(t) \neq 0$, since it is the Alexander  
               polynomial $\AlxK(t)$ of the knot $K$.            
        \phcc  
               $\AlxL(t) = \detb \rt  M - \rti M^T & \rtt V^T \\ 
                                 \rtt V            & \rtt A   \\ \enddet$
        \phcd          $ = \detb R(t) & \rtt V^T              \\ 
                                  0   & \rtt A - \rtt^2 \VRV  \\ \enddet $
        \phcd          $ = \rtt^m \ \detb R(t) & V^T            \\ 
                                         0   & A - \rtt \VRV  \\ \enddet $   
        \phcd          $ = \rtt^m \ \det R(t) \ \det (A - \rtt \VRV) $. 
        \pvcn  Thus we have
               $$\AlxL(t)= \rtt^m \AlxK(t) \det(A  - \rtt \VRV).$$ 
        \pvcn  It follows from \cite[Theorem 4.1]{PY} that 
               $$\lkx(\JipF,\JjF)-\lk(\Jip, \Jj)=
                (1-t) \VJi (tM-M^T)^{-1} \VJj^T.$$
        \pvcn  Since $R(t) = \rti(t M - M^T)$, we have that 
               $$R^{-1}(t) = \rt (t M - M^T)^{-1}.$$
        \pvcn  Therefore the $(i,j)$-entry of the matrix $A - \rtt \VRV$ is
        \pvcn  $$\si(r_i,r_j) - \rtt \VJi R^{-1}(t) \VJj^T
                =\lk(\Jip,\Jj)+(1-t)\VJi (tM-M^T)^{-1} \VJj^T
                =\lkx(\JipF,\JjF).$$
                \begin{flushright}$\Box$\end{flushright}

%
%
%
%
%
%
%
%
%
%
 
                                 \section{Calculation of $\bekl$}
                                 \label{sec:calculation}

%
%
%
%
%
%
%
%
%

               In \cite{TC-CMH85}, an algebraic method of computing 
               $\lk(J^{F(k)+},J^{F(k)})$ is given. In a similar way we can
               calculate $\bekl$. Also as its corollaries we have that
               $\bekl=(-1)^l \be_F^{k+l,0}$ (Corollary \ref{cor:beta}) and
               the Taylor expansion of the inverse of the Alexander
               matrix around $t=1$ (Corollary \ref{cor:Alexander}).

\begin{Thm}   $($c.f. \cite[\S 8]{TC-CMH85}$)$ 
               \label{thm:betabymatrices}        
               Let $\KJaJb$ a $3$-component link with $\lk(K,\Ji) =0$ $\iab$ 
               and let $F$ be a Seifert surface for $K$ in $\cmpl (\JaJb)$.
               Let $\om = \pq$ be a symplectic basis for $H_1(F)$ and $M$ 
               the Seifert matrix with respect to $\om$. Then we have the 
               following, where $\VJi=(\jipa, \cds, \jipg, \jiqa, \cds, \jiqg)$
               and $-P$ is the intersection matrix $M-M^T$ with respect to $\om$: 
               $$\bekl(\Ja,\Jb) =   
               (-1)^l \VJa (PM)^{k+l-1} P \VJbt \ (k \geq 1, l \geq 0).$$ \end{Thm}
          \bgp 
               Let $u^{\pm}$ be $\mtxf \pa^\pm & \cds & \pg^\pm & \qa^\pm 
               & \cds & \qg^\pm \\ \endmtx^T$. Then from the proof of 
               Proposition \ref{prop:derivative} we have:
         \pvcn 
               $\pJia   = \sum_{j=1}^g (\jiqj \pjp   - \jipj \qjp)$
         \phc         $ = \mtxf \jipa & \cds & \jipg & 
                                \jiqa & \cds & \jiqg \\ \endmtx
                          \mtxf -\qap & \cds & -\qgp & 
                                 \pap & \cds & \pgp  \\ \endmtx^T $   
         \phc         $ = \VJi  \mtxb   0   & -I_g \\
                                      I_g  &  0   \\ \endmtx
                               \mtxf  \pap & \cds & \pgp  & 
                                      \qap & \cds & \qgp  \\ \endmtx^T$
         \phc         $ = \VJi P u^+$,
         \pvcn and
         \pvcn $V_\pJia = \mtxf \lk(\VJi P u^+, \pa) & \cds & 
                                \lk(\VJi P u^+, \pg) &
                                \lk(\VJi P u^+, \qa) & \cds & 
                                \lk(\VJi P u^+, \qg) \\ \endmtx$
         \phc         $ = \VJi P \mtxf 
                                \papa & \cds & \papg & \paqa & \cds & \paqg \\
                                \vds  &      & \vds  & \vds  &      & \vds  \\
                                \pgpa & \cds & \pgpg & \pgqa & \cds & \pgqg \\
                                \qapa & \cds & \qapg & \qaqa & \cds & \qaqg \\
                                \vds  &      & \vds  & \vds  &      & \vds  \\
                                \qgpa & \cds & \qgpg & \qgqa & \cds & \qgqg \\ 
                                \endmtx$
         \phc         $ = \VJi P M.$
         \pvcn In the same way we obtain that $\pJin = V_\pJim P u^+$ and
               $V_\pJin = V_\pJim P M$.  Therefore inductively we obtain that   
                 $\pJin = \VJi (PM)^{n-1} P u^+$. Similarly we also have that
                 $\nJin = \VJi (PM^T)^{n-1} P u^-$. 
               Therefore we can calculate $\beaz$ and 
               $\bekl$ ($k,l \geq 1$) as follows:
        \pvcn  $\beaz(\Ja,\Jb) = \lk(\Ja^{(1)+}, \Jb) = \lk(\VJa Pu^+,\Jb)$
        \phd                $  = \VJa P(\lk(p_1^+,\Jb), \cds, \lk(p_g^+,\Jb),
                                        \lk(q_1^+,\Jb), \cds, \lk(q_g^+,\Jb))$
        \phd                $  = \VJa P \VJb^T.$
        \pvcn  and
        \pvcn  $\bklab = \lk(\pJak, \nJbl)
                       = \lk(\VJa (PM)^{k-1} P u^+, \VJb (PM^T)^{l-1} P u^-)$
              \phcdh $ =    (\VJa (PM)^{k-1} P)  M (\VJb (PM^T)^{l-1} P)^T  $
              \phcdh $ =     \VJa (PM)^{k-1} P   M P^T (M P^T)^{l-1} \VJb^T $   
              \phcdh $ =     \VJa (PM)^k  (P^T M)^{l-1} P^T \VJb^T.$   
        \pvcn Now since $P^T=-P$, we complete the proof.
              \end{proof}

\begin{Cor}   \label{cor:beta} We have that
              $\bekl(\Ja,\Jb) = (-1)^l \be_F^{k+l,0}(\Ja,\Jb).$  \end{Cor}

\begin{Cor}   \label{cor:Alexander}
              Let $M$ be a Seifert matrix for a knot with respect to 
              a symplectic basis. Then the inverse of the Alexander matrix 
              $tM-M^T$ may have the following expansion:
              $$(tM-M^T)^{-1} = \sumnzi (1-t)^n (PM)^n P.$$ \end{Cor}
          \bgp 
              Let $F$ be a Seifert surface for the knot and let $\JaJb$ be
              a link in $\cmpl F$. By combining \cite[Theorem 4.1]{PY} and 
              Theorem \ref{thm:ICC}, we have that
              $$\lkxabf = \lk(\Ja,\Jb) + (1-t) \VJa (tM - M^T)^{-1} \VJbt
                        = \sumkzi (-1)^k \alfk (t-1)^k.$$
\noindent     Since $\lk(\Ja,\Jb)=\azab$, we have
              $$\VJa (tM - M^T)^{-1} \VJbt = \sumkai (-1)^{k-1} \alfk (t-1)^{k-1}.$$
              \pvcn 
              From Theorem \ref{thm:betabymatrices},
              $$\alfn = \VJa (PM)^{n-1} P \VJbt \ (n \geq 1). $$
\noindent     Hence we have
              \pvcn
              $(tM - M^T)^{-1} = \sumkai (-1)^{k-1} (t-1)^{k-1} (PM)^{k-1} P$
        \phcb $= \sumnzi (-1)^n (t-1)^n (PM)^n P.$
              \end{proof}

%
%
%
%
%
%
%
%
%
%
%
%
%
%
%



\begin{thebibliography}{99}

          \bibitem{TC-CMH85}        T.D. Cochran 
                               {\em Geometric invariants of link cobordism}, 
                                    Comment. Math. Helvetici 
                               {\bf 60} (1985), 291--311.
          \bibitem{TC-IM85}         T.D. Cochran 
                               {\em Concordance invariance of coefficients 
                                    of Conway's link polynomial}, 
                                    Invent. Math. 
                               {\bf 82} (1985), 527--541.
          \bibitem{FH-OMJ58}        F. Hosokawa   
                               {\em On $\na$-polynomial of links}, 
                                    Osaka Math. J. 
                               {\bf 10} (1958), 273--282.
          \bibitem{JH-PAMS85}       J. Hoste    
                               {\em The first coefficient of the Conway polynomial}, 
                                    Proc. AMS. 
                               {\bf 95} (1985), 299--302.
          \bibitem{JL-Knots96}      J.P. Levine  
                               {\em The Conway polynomial of 
                                    an algebraically split link}, 
                                    Proceedings of Knots 96, ed. S. Suzuki, 
                                    World Sci. Puli. Co. (1997), 91--98. 
          \bibitem{JL-CMH99}        J.P. Levine 
                               {\em A factorization of the Conway polynomial}, 
                                    Comment. Math. Helvetici 
                               {\bf 74} (1999), 27--53.
          \bibitem{KY-IM79}         S. Kojima and M. Yamasaki 
                               {\em Some new invariants of links}, 
                                    Invent. Math. 
                               {\bf 54} (1979), 213--228.
          \bibitem{PY}              J.H. Przytycki and A. Yasuhara 
                               {\em Linking numbers in rational homology $3$-spheres, 
                                    cyclic branched covers and infinite cyclic covers}, 
                                    Trans. AMS.
                               {\bf 356} (2004), 3669--3685.
          \bibitem{NS-TA84}         N. Sato     
                               {\em Cobordisms of semi-boundary links},
                                    Topology and its Applications 
                               {\bf 18} (1984), 225--234.

                                  \end{thebibliography}
                                      \end{document}